\documentclass{article}
\usepackage[a4paper, total={6in, 10in}]{geometry}
\usepackage{amsmath,amsthm,amssymb,mathrsfs,amstext, titlesec,enumitem, comment, graphicx, color, xcolor, inputenc, comment}
\usepackage{hyperref}
\usepackage{xpatch}
\makeatletter
\def\Ddots{\mathinner{\mkern1mu\raise\p@
\vbox{\kern7\p@\hbox{.}}\mkern2mu
\raise4\p@\hbox{.}\mkern2mu\raise7\p@\hbox{.}\mkern1mu}}
\makeatother

\titleformat{\section}
{\normalfont\Large\bfseries}{\thesection}{1em}{}
\titleformat*{\subsection}{\Large\bfseries}
\titleformat*{\subsubsection}{\large\bfseries}
\titleformat*{\paragraph}{\large\bfseries}
\titleformat*{\subparagraph}{\large\bfseries}

\makeatother
\newtheorem{theorem}{Theorem}[section]

\theoremstyle{definition}

\newcommand{\bN}{{\mathbb N}}
\newcommand{\bZ}{{\mathbb Z}}

\date{\vspace{-5ex}}

\begin{document}

\title{ Monochromatic Polynomial sumset structures on $\bN$: an ultrafilter proof}
\author{Sayan Goswami \footnote{Ramakrishna Mission Vivekananda Educational and Research Institute, Belur Math,
Howrah, West Benagal-711202, India, email: \textit{ sayan92m@gmail.com}.}}

\maketitle
\begin{abstract}

Recently in \cite{LX}, using machinery's from Ergodic theory,  Z. Lian, and R. Xiao proved if $P$ is any polynomial with no constant term, then for every finite coloring of $\bN$, there exists two infinite subsets $B,C$ of $\bN$ such that the set $\{P(b)+P(c):b\in B, c\in C\}$ is monochromatic. In this article we improve their result by proving that  instead of taking such polynomials we can choose any function $f$ having the property that $f(\bN)\setminus \bN$ is finite. We use ultrafilter techniques to prove our result.
\end{abstract}
Keywords: Polynomial sumset, Algebra of the Stone-\v{C}ech compactification of discrete semigroups\\
Subject class: 05D10, 54D35.

\section{Introduction}

A core problem in arithmetic Ramsey
theory is the characterization of families $\mathcal{F}$ of
subsets of a semigroup $(S,\cdot )$ that are partition regular, $\textit{i.e.}$
the families have the property that whenever $S=\bigcup_{i=1}^{r}A_{i}$
is a finite coloring of $S$, at least one of the $A_{i}$ is in $\mathcal{F}$. For any $r$-coloring of $S,$ 
$A\subset S$ will be called monochromatic if each member of $A$ is in the same
color. A cornerstone result in arithmetic Ramsey theory is van
der Waerden's theorem \cite{2}, which states that for any $l,r\in\mathbb{N}$,
and for any $r$ coloring of $\mathbb{N},$ there exists a monochromatic 
arithmetic progression (AP) of length $l$. Polynomial extensions of classical Ramsey theoretic results are much difficult to prove. In \cite{1},  V. Bergelson and A. Leibman found the polynomial extension of this theorem. Later in \cite{M}, J. Moreira, F. Richter, and D. Robertson proved the frollowing theorem which was a conjecture of Erd\"{o}s.
 
\begin{theorem}
    If $A\subseteq \bN$ has positive upper density; i.e. $\limsup_{N\rightarrow\infty}\frac{|A\cap |1,2,\ldots ,N|}{N}>0$, then $A$ contains $B+C=\{b+c:b\in B,c\in C\}$ for some infinite sets $B,C\subseteq \bN.$
\end{theorem}

Recently in \cite{LX},  Z. Lian, and R. Xiao studied a polynomial variant of the above problem. Using Ergodic theory, they proved the following result.

\begin{theorem}\cite[Theorem 1.6]{LX}\label{1}
  For any finite coloring of $\bN$, and any polynomial with positive leading coefficient and zero constant term,  there exist infinite subsets $B,C\subseteq \bN$ such that $B\subseteq P(\bN)=\{P(n):n\in \bN\}$, and $B+P(C)=\{b+P(c):b\in B, c\in C\}$ is monochromatic.
\end{theorem}

This article is aimed to improve Theorem \ref{1}. Using the theory of ultrafilters we prove the following theorem, which is more general.

\begin{theorem}\label{2}
    Let $f:\bZ \rightarrow \bZ$ be any function such that $f(\bN)\setminus \bN$ is finite. Then for every finite partition of $\bN$, there exists two infinite subsets $B,C\subseteq \bN$ such that $f(B)+f(C)=\{f(b)+f(c):b\in B, c\in C\}$ is monochromatic.
\end{theorem}

\subsection{Preliminaries}

In this section we recall some basic prerequisites from the ultrafilter theory. For details we refer the book \cite{HS} to the readers. 

A filter $\mathcal{F}$ over any nonempty set $X$ is a collection of subsets of $X$ such that

\begin{enumerate}
    \item $\emptyset \notin \mathcal{F}$, and $X\in \mathcal{F}$,
    \item $A\in \mathcal{F}$, and $A\subseteq B$ implies $B\in \mathcal{F},$
    \item $A,B\in \mathcal{F}$ implies $A\cap B\in \mathcal{F}.$
\end{enumerate}
Using Zorn's lemma we can guarantee the existence of maximal filters which are called ultrafilters. Any ultrafilter $p$ has the following property:
\begin{itemize}
    \item if $\bN=\cup_{i=1}^rA_i$ is any finite partition of $\bN$, then there exists $i\in \{1,2,\ldots ,r\}$ such that $A_i\in p.$
\end{itemize}
Let  $(S,\cdot)$ be any discrete semigroup. Let $\beta S$ be the collection of all ultrafilters. For every $A\subseteq S,$ define $\overline{A}=\{p\in \beta S: A\in p\}.$ Now one can check that the collection $\{\overline{A}: A\subseteq S\}$ forms a basis for a topology. This basis generates a topology over $\beta S.$  We can extend the operation $``\cdot ''$ of $S$ over $\beta S$  as: for any $p,q\in \beta S,$ $A\in p\cdot q$ if and only if $\{x:x^{-1}A\in q\}\in p.$ With this operation $``\cdot "$, $(\beta S,\cdot)$ becomes a compact Hausdorff right topological semigroup. One can show that $\beta S$ is nothing but the Stone-\v{C}ech compactification of $S.$ If $S,T$ are two discrete semigroups, then Now one can easily extend any function from $S$ to $T$ to the function from $\beta S$ to $\beta T$ as follows:
\begin{itemize}
    \item if $f:S\rightarrow T$ is a function then the function  $\tilde{f}:\beta S\rightarrow \beta T$ is defined as $A\in \tilde{f}(p)$ if and only if $f^{-1}(A)=\{x\in S: f(x)\in A\}\in p.$
\end{itemize}

\section{Proof of our Theorem \ref{2}}

\begin{proof}[\textbf{Proof of Theorem \ref{2}:}]
Let $f:\bZ \rightarrow \bZ$ be any function such that it $f(\bN)\setminus \bN$ is finite. We claim that for every nonprinciple ultrafilter $p\in \beta \bN$, we have  $\tilde{f}(p)\in \beta \bN$, and it is a nonprinciple ultrafilter in $ \beta \bN$. To prove our claim let $\tilde{f}:\beta \bZ\rightarrow \beta \bZ$ be the continous extension of $f$. Passing to the restriction of $\beta \bN$, we may think that $\tilde{f}$ is a map from $\beta \bN$ to $\beta \bZ$. Let $p\in \beta \bN$ be nonprinciple. So $\tilde{f}(p)$ is nonprinciple in $\beta \bZ$. Now we need to show that $\bN\in \tilde{f}(p)$. Note that $f(\bN)\in \tilde{f}(p).$ But $f(\bN)=(f(\bN)\setminus \bN)\cup (f(\bN)\cap \bN).$ As $f(\bN)\setminus \bN$ is finite, and $\tilde{f}(p)$ is nonprinciple, we have $f(\bN)\cap \bN\in \tilde{f}(p).$ That is $\bN\in \tilde{f}(p).$ This completes the proof of the claim.

Let $A\in \tilde{f}(p)+\tilde{f}(p).$ Then $B=\{r: -r+A \in \tilde{f}(p)\}\in \tilde{f}(p).$ Choose $f(x_1)\in B.$ Then $-f(x_1)+A\in \tilde{f}(p).$ Hence $ B\cap (-f(x_1)+A)\in \tilde{f}(p).$ Now choose $f(y_1)\in B\cap (-f(x_1)+A).$ Then $f(x_1)+f(y_1)\in A.$ Now for some $N>1$, assume we have elements $\langle x_n\rangle_{n=1}^N$, and $\langle y_n\rangle_{n=1}^N$ such that:

\begin{enumerate}
    \item $-f(x_i)+A\in \tilde{f}(p)$, and $-f(y_j)+A\in \tilde{f}(p)$ for every $i,j\in \{1,2,\ldots ,N\}$,
    \item for every $i,j\in \{1,2,\ldots ,N\}$, we have $f(x_i)\in \cap_{j=1}^N (-f(y_j)+A)$, and $f(y_i)\in \cap_{j=1}^N (-f(x_j)+A)$.
\end{enumerate}

As for every $i\in \{1,2,\ldots ,N\}$, $-f(x_i)+A\in \tilde{f}(p)$, we have $ \bigcap_{i=1}^N -f(x_i)+A\in \tilde{f}(p)$. Now choose $f(y_{N+1})\in B\cap \bigcap_{i=1}^N (-f(x_i)+A).$ Again from the definition of $B,$ we have $-f(y_{N+1})+A\in  \tilde{f}(p).$ Again $\bigcap_{i=1}^{N+1} (-f(y_i)+A)\in \tilde{f}(p).$ Now choose $f(x_{N+1})\in B\cap \bigcap_{i=1}^{N+1} (-f(y_i)+A).$
Clearly our new $x_{N+1},$ and $y_{N+1}$ satisfies induction hypothesis.
So we have two sequences $\langle x_n\rangle_{n\in \bN}$, and  $\langle y_n\rangle_{n\in \bN}$ such that for every $N\in \bN$, and $i,j\in \{1,2,\ldots ,N\}$, we have $f(x_i)\in \cap_{j=1}^N (-f(y_j)+A)$, and $f(y_i)\in \cap_{j=1}^N (-f(x_j)+A)$. Now assume $s,t\in N$. Choose $N=\max \{s,t\}.$ Then clearly $s,t\in \{1,2,\ldots ,N\},$ and so we have $f(x_s)+f(y_t)\in A.$ 

This completes the proof.

\end{proof}

\section*{Acknowledgement} The author is supported by NBHM postdoctoral fellowship with reference no:\\ 0204/27/(27)/2023/R \& D-II/11927.

\end{document}